\documentclass[10pt]{amsart}
\usepackage[cp1251]{inputenc}
\usepackage[english]{babel}
\usepackage{amsmath}
\usepackage{amssymb}
\usepackage{amsfonts}

\begin{document}

\thispagestyle{empty}
\title{}

\begin{center}

\end{center}
\maketitle

\vskip-6mm

K. V. Storozhuk\vskip2mm

\centerline{\bf Symmetric Invariant Subspaces} \centerline{\bf of
Complexifications of Linear Operators\rm}\vskip1mm

\bf 1. Introduction.\rm\vskip1mm

Let $T:X\to X$ be a linear operator on a complex Banach space.
Denote by $\sigma(T)$ and $R(\lambda, T)$ its spectrum and
resolvent.

Suppose that $\sigma(T)$ is disconnected, $\sigma(T)=F_1\cup F_2,$ where
$F_1$ and $F_2=\sigma(T)\backslash F_1$ are some components of
spectrum. Let $\gamma$ be a contour surrounding $F_1$. Consider the
spectral projection $P=\frac{1}{2\pi i}\int_\gamma R(\lambda,
T)d\lambda $. Its image $[F_1]$ and kernel $[F_2]$ are both invariant
subspaces and $\sigma (T|_{[F_i]})=F_i,\ i=1,2$.

Suppose that the spectrum is connected. Cutting out the subset
$F\subset \sigma(T)$ by a contour $\gamma$, we must multiply
the resolvent to weight function $g$ that is small in some neighborhood
of the points of intersection $\gamma\cap F$: $f(T)=\frac{1}{2\pi
i}\int_\gamma R(\lambda, T)g(\lambda)d\lambda$.
This method allows to obtain analogs of spectral subspaces with
constraints on the growth of the resolvent near the spectrum (see \cite{Da,L}).
These constraints are fulfilled, for example, if the powers of
$T^{\pm n}$ increase slowly. If the spectrum lies on the unit sircle
$\Lambda$ then the condition of non-quasianalyticity
$\frac{\sum_{n=-\infty}^\infty\|T^n\|}{1+n^2}<\infty$ is equivalent
to the Levinson condition, an integral constraint on the growth
of the resolvent near the spectrum; such operators also have
a separated spectrum \cite{LMF}.

Suppose now that $T:X\to X$ and $X$ is real.  The spectral
projection corresponding to the symmetric component of the spectrum
of the complexification $T_\Bbb C:X_\Bbb C\to X_\Bbb C$ gives us a
symmetric invariant subspace $L_\Bbb C$. Its real part $L\subset X$
is $T$-invariant. See, for example, theorem 5.3 of \cite{BZ}.

Even in the case where the spectrum of $T_\Bbb C$ is connected, it
is still easy to get a symmetric $T_\Bbb C$-subspace by using the
above method. We must integrate over a \it symmetric \rm contour
(with respect to the real axis) with a \it symmetric \rm  function
$g$. The ``real part'' of $f(T_\Bbb C)$ is a $T$-invariant subspace
in $X$. In the second part of this note, we will expose a
``realification'' of one of such methods. As an application, we
receive a theorem whose complex version was proved in \cite{We}. The
last statement of the ``complex'' theorem is used quite often and
goes back to Theorem~J of \cite{Go}.\vskip1mm

\bf Theorem 1\it.  Suppose that $T:X\to X$, $X$ is a real Banach
space, $T$ is an invertible linear operator such that
$\|T^n\|=o(|n|^k)$ as $n\to\pm\infty$. If $\dim X>2$ then $T$ has an
invariant subspace. In particular, a linear isometry $T:X\to X$ of a
real space has an invariant subspace if $\dim\,X>2$\rm.\vskip1mm

Perhaps these results are already known, but the author has not
found any direct reference. Let us observe, however, that theorems
about invariant spaces of compact operators were generalized to the
real case. See, for example, \cite{AAST} and the references therein.

Of course, not all invariant subspaces can be obtained by spectral
methods. There are operators of Voltrerra type such that, for each
invariant subspace $L$, $\sigma(T|_L)=\sigma(T)=[0,1]$  (see
\cite{LM1}). Generally speaking, if the operator $T_\Bbb C$ has
invariant subspaces it is not known whether there are symmetric
subspaces among them; the existence of such subspaces is equivalent
to the validity of Conjecture 3 in \cite{AAST}.\vskip2mm

\bf 2. Detailed Definitions and the Proof of the Theorem
1\rm.\vskip2mm

Below by an \it invariant subspace \rm we mean a \it
nontrivial closed \rm invariant subspace taken to itself by an operator $T:X\to X$.

Suppose that $T:X\to X$, $X$ is a complex Banach space, and $x\in
X$. The map $\lambda\mapsto R(\lambda, T)x$ is an $X$-valued
function holomorfic outside $\sigma(T)$. If this map has a
single-valued analytic extension to some set $\rho(x)$ then the set
$\sigma(x):=\Bbb C\backslash \rho(x)\subset \sigma(T)$ is called the
local spectrum of $x$ and the corresponding extension is called the
local resolvent of $x$.

Now, suppose that $X$ is real  space. The complexification of $X$ is
the space $X_\Bbb C$ of elements of the form $(x+iy)$; it is natural
to call the vectors $x,y\in X$ \,\, real $(\mbox{Re}(z))$ and
imaginary $(\mbox{Im}(z))$ parts of $z$. The space $X_\Bbb C$ is
endowed by the conjugation $J:X_\Bbb C\to X_\Bbb C$,
$J(x+iy)=(x-iy)$. Introduce a norm by $\|z\|^2=\max\{\|Re(\lambda
z)\|^2+\|Im(\lambda z)\|^2\mid \lambda\in\Bbb C,\ |\lambda|=1\}.$
This norm is equivalent to the norm of the direct sum  $X\oplus X$.

The complexification $T_\Bbb C:X_\Bbb C\to X_\Bbb C$ of an operator
$T:X\to X$ is defined by $T_\Bbb C(x+iy)=(Tx+iT_y)$.

We call a subset $F\subset \Bbb C$ {\it symmetric} if $F$ is
symmetric with respect to the real line, i.e. $F=\overline{F}$.
Similarly, a subset $Z\subset X_\Bbb C$ is symmetric if $J(Z)=Z$. It
easy to see that if $Z\subset X_\Bbb C$ is a symmetric subspace then
$Z=L_\Bbb C$, where $L=\mbox{Re}\,Z=\mbox{Im}\,Z\subset X$.

\bf Lemma\it. Suppose that an operator $ T_\Bbb C:X_\Bbb C\to X_\Bbb
C$ is the complexification of some real operator $T:X\to X$. Then
the spectrum $\sigma(T_\Bbb C)\subset\Bbb C$ is symmetric. If
$T_\Bbb C$ admits a local resolvent, then
$\sigma(J(z))=\overline{\sigma(z)}$ for each $z\in X_\Bbb C$.
\rm

Proof\rm. It is easy to verify the equality $R(\bar{\lambda}, T_\Bbb
C)=J\circ R(\lambda,T_\Bbb C)\circ J$, and so the spectrum of
$T_\Bbb C$ is symmetric (this is lemma 4.1 of \cite{BZ}). If $z\in
X_\Bbb C$ and a function $f$ is an analytic extension of the
resolvent $R(\lambda, T)z$ then the function $J\circ f\circ J:
\lambda\to J\circ R(\lambda,T)(J(x))$ is an analytic extesion of the
resolvent $R(\bar{\lambda}, T_\Bbb C) x$. Therefore, both maximal
extensions coincide and $\sigma(J(z))=\overline{\sigma(z)}$. The
lemma is proved.

Now we prove Theorem 1. Let $F$ be a \it symmetric \rm arc
containing a part of the spectrum of $T_\Bbb C$ and let $[F]\subset
X_\Bbb C$ be the space of vectors whose local spectrum is included
in $F$. The spectrum of $T_\Bbb C$ is separated \cite{LMF} and hence
the subspace $[F]$ is $T_\Bbb C$-invariant. Lemma 1 implies that
this subspace is symmetric. As is easy to verify,
$\mbox{Re}[F]\subset X$ is a $T$-invariant subspace.



The spectrum $T_\Bbb C$ may consist of at most two points
$\eta,\bar{\eta}\in\Lambda$; therefore, no symmetric arc $F$
contains a \it part \rm of the spectrum. In this situation, the
separatedness of the operator and the constraint on the growth of
$\|T^{\pm n}\|$ allow us to use the Gelfand---Hille theorem
involving which it is easy to conclude that $((T_{\Bbb
C}-\eta)(T_{\Bbb C}-\bar{\eta}))^{k+1}=0.$ Hence, the closure of the
range of the operator $(T_{\Bbb C}-\eta)(T_{\Bbb C}-\bar{\eta})$ is
not equal to $X_\Bbb C$. The coefficients of the polynomial $T_{\Bbb
C}^2-aT_{\Bbb C}+bI=(T_{\Bbb C}-\eta)(T_{\Bbb C}-\bar{\eta})$ are
real, and so the closure of the range $(T^2-aT+bI)$ is not equal to
$X$. Cf. the proof of Theorem~3 in [5]. Therefore, either it is an
invariant subspace or $(T^2-aT+bI)X=0$. In the latter case, every
vector $x\in X$ generates at most a two-dimentional invariant
subspace.

Consider the case where $T$ is an isometry. If $T$ is bijective then $\|T^{\pm
n}\|=1$ for all $n\in \Bbb N$ and everything is proved. If $TX\neq X$
then $TX$ is a desired invariant subspace. Theorem 1 is proved.

{\bf Remark 1.} We use very general results of \cite{LMF}, but for
proving the separatedness of $[F]$ in our theorem, it would have
sufficed to refer to Leaf's work \cite{L} and even Dunford's basic
result (see Corollary 9 of Chapter XVI \S 5 \cite{D-Sh}).

{\bf Remark 2.} At the end of article \cite{LMF}, it is mentioned
that Wermer has shown the existence of invariant subspaces, ``of course,
if the spectrum contains more than one point''. This is not the case;
in \cite{We}, the one-point case is considered separately in the proof of
theorem 3  in \cite{We}. That is the argument we have used in the proof of
the Theorem 1 in the case of a two-point spectrum.

{}
\end{document}